\documentclass[12pt]{amsart}
\usepackage[totalwidth=14cm,totalheight=21cm]{geometry}
\newtheorem{remark}{Remark}[section]

\newtheorem{theorem}{Theorem}[section]

\begin{document}
\title{On stable  constant mean curvature surfaces in $\mathbb S^2\times \mathbb R$ and $\mathbb H^2\times \mathbb R $}

\date{April  3rd, 2008, (v2)}

\author{Rabah SOUAM}
\address{Institut de Math{\'e}matiques de Jussieu,  CNRS  UMR 7586, Universit{\'e} Paris Diderot - Paris 7,
"G{\'e}om{\'e}trie et Dynamique ", Site Chevaleret,  Case 7012,  75205- Paris Cedex 13, France.}
\email{souam@math.jussieu.fr}

\newtheorem{proposition}{Proposition}[section]

\begin{abstract}
 We study the stability of immersed compact constant mean curvature (CMC)
 surfaces without boundary in some Riemannian 3-manifolds, especially  in the Riemannian product spaces  $\mathbb S^2 \times \mathbb R$ and $\mathbb H^2\times\mathbb
R.$  We prove that rotational CMC spheres in $\mathbb H^2\times\mathbb R$ are all stable, whereas in $\mathbb S^2\times\mathbb R$ there exists some
value $H_0\approx 0.18$ such that rotational CMC spheres are stable for $H\geq H_0,$  and unstable for $0<H<H_0.$ We show that a compact
stable immersed CMC surface in $\mathbb S^2\times \mathbb R$ 
 is either a finite union of horizontal slices or a rotational sphere. In the more general case of  an ambient manifold
 which is a simply connected conformally flat three manifold with nonnegative  Ricci curvature  we show that a closed stable
immersed CMC surface is either a sphere or  an embedded torus.  Under the weaker  assumption that the  scalar  curvature is nonnegative, we prove that a closed stable immersed CMC surface
 has genus at most three. In  the case of $\mathbb H^2\times \mathbb R$
we show that a closed stable immersed CMC surface is a rotational sphere if it has mean curvature  $H\ge 1/\sqrt {2}$
and that it has genus at most one if $1/\sqrt{3} < H < 1/\sqrt {2}$  and 
 genus at most two if $H=1/\sqrt{3}.$
\end{abstract}

\maketitle

\section{INTRODUCTION AND PRELIMINARIES}
Compact surfaces of constant mean curvature (CMC) in  Riemannian 3-manifolds  are critical points of a variational problem:
 minimizing the area
under a volume constraint. It is therefore natural to study the {\it stable} ones, that is those minimizing the area functional
 up to second order for {\it volume-preserving} deformations \cite{barbosa et al}.  Let $\Sigma$ be an  immersed compact 
surface without boundary and  with constant mean curvature $H$  in a  Riemannian
3-manifold $(M,\langle,\rangle).$ Suppose $\Sigma$ is  two-sided, i.e there is a global unit normal field $N$ along it. 
This condition is always satisfied when $H\neq 0,$ and is equivalent to the orientability of $\Sigma$ when the ambient manifold $M$ is
orientable. This will be the case in this article since we will consider simply connected ambient manifolds. Also we can always choose the
orientation on the surface so that its (constant) mean curvature is nonnegative, therefore we can always suppose $H\geq 0.$

We consider on $\Sigma$ the following quadratic form:
$$ Q(u,v)= \int_{\Sigma}\langle \nabla u, \nabla v\rangle - (|\sigma|^2+\text{Ric}(N))uv \,dA=-\int_{\Sigma}
 u Lv\,dA,\quad u,v\in\mathcal{C}^\infty(\Sigma),$$
  where $\nabla u$ stands for the gradient of $u,$   $\sigma$ is the second fundamental form of 
the immersion and Ric$(N)$ denotes the Ricci curvature of $M$ evaluated on the  field $N.$  The linear operator 
$L=\Delta + |\sigma|^2+\text{Ric}(N)$ is the Jacobi operator of the surface, $\Delta$ being the Laplacian on $\Sigma.$
 The stability
condition means that: 
\begin{equation}\label{stability}
Q(u,u)\geq 0,\quad \text{for any}\, u\in  \mathcal  {C}^\infty(\Sigma) \quad \text{satisfying} \int_{\Sigma} u\,dA =0.
\end{equation}
 See  \cite{barbosa et al} for the details. 

The Jacobi operator $L$ is also the linearized operator of the mean curvature functional. This means that if 
$X:\Sigma\longrightarrow (M, \langle,\rangle)$ denotes a CMC immersion as above, then for any smooth deformation $X_t$ 
through immersions of $X=X_0,$ the derivative of the mean curvature functional at $t=0$ is given by:

\begin{equation}\label{linearization} 
2\frac{\partial H_t}{\partial t}|_{t=0} = L (\langle \frac{\partial X_t}{\partial t}|_{t=0}, N\rangle).
\end{equation}

Stability of CMC surfaces is an important notion in solving the {\it isoperimetric problem} in Riemannian manifolds 
as the boundary of
 an {\it isoperimetric region} is a stable CMC surface. Recent advances 
on the isoperimetric problem can be found in the paper by Ros
\cite{ros1}. 

Stability of CMC surfaces has attracted much attention. The first results in this direction are due to Barbosa,  do Carmo
\cite{barbosa-do carmo}, and  Barbosa,  do Carmo,  Eschenburg \cite{barbosa et al}, who established that in simply connected space forms, 
the
round spheres are the only stable immersed compact CMC surfaces, a result which holds in all dimensions. In this paper, we will be
concerned with  stability of compact surfaces in some Riemannian 3-manifolds, mainly the Riemannian product 
 spaces $\mathbb S^2\times\mathbb R$ and $\mathbb H^2\times\mathbb R.$ The theory of CMC surfaces in these spaces is currently an 
 active research field. Abresch and  Rosenberg \cite{abresch-rosenberg} have proved that an immersed CMC sphere in
 any of these spaces is 
 necessarily rotationally invariant.  The rotational CMC spheres are embedded and form, in both spaces, 
a 1-parameter family-up to ambient isometries- parameterized by their mean curvature. They have first  appeared in the study of the
isoperimetric problem, 
 by Hsiang and Hsiang 
 \cite{hsiang} in $\mathbb H^2\times\mathbb R$ and by  Pedrosa in $\mathbb S^2\times\mathbb R$ \cite{pedrosa}.  In $\mathbb H^2\times\mathbb R,$ 
the solutions to the isoperimetric problem are precisely the domains bounded by the rotational CMC spheres. In $\mathbb S^2\times\mathbb R,$ 
Pedrosa has shown that the rotational CMC spheres bound isoperimetric domains only 
  when they have mean curvature 
 $H\geq H_1\approx 0.33.$
 
 In section 1, we establish our first result. It concerns the stability of rotational CMC spheres (Theorem \ref{stablerotational}). In $\mathbb H^2\times\mathbb R,$
 we prove they are all stable. In $\mathbb S^2\times\mathbb R,$
 our study reveals a rather unexpected phenomenon: we show that rotational CMC spheres are stable if and only if they have 
 mean curvature $H\geq H_0\approx 0.18.$ The value $H_0$ is the (unique) one where the area of these surfaces- also the
  enclosed volume- reaches its maximum and these spheres are stable (resp. unstable)
 as long as the
 area- also the enclosed volume- is decreasing (resp. increasing) as a function of the mean curvature. In particular, the rotational CMC
spheres are still stable when their mean curvature lies in the interval $[ H_0,  H_1[ $ although they do not bound isoperimetric domains. 

In Section 2, we obtain some restrictions on the topology of stable CMC surfaces in some class or Riemannian 3-manifolds
under  positive curvature assumptions. Improving previous results by several authors, Ros \cite{ros2} proved  that an
 immersed closed  stable two-sided CMC surface in a 3-manifold with nonnegative Ricci curvature has genus at most 3. The result
is sharp and cannot be improved in general since  Ross \cite{ross} showed that the classical Schwarz-$P$ surface, of genus 3 in the
flat cubic lattice, is stable. However one can hope to improve this bound in some classes of manifolds. Recently,  Ros \cite{ros3}
considered the case of a product manifold $M\times\mathbb R$ where $M$ is  a complete, orientable surface with nonnegative
Gaussian curvature and proved that if $\Sigma$ is a compact orientable stable and {\it embedded} CMC surface in $M\times\mathbb R,$
then either $M$ is compact and $\Sigma$ a finite union of horizontal slices, or $\Sigma$ is a connected surface with genus $\leq 2.$
We consider here the case of an  ambient manifold which is simply connected and conformally flat with nonnegative 
Ricci curvature. We prove (Theorem \ref{thm:conformal}) that a compact orientable stable CMC in such a manifold 
 is either a finite union of totally geodesic surfaces, a sphere or an embedded torus. In particular, if we assume the Ricci curvature is positive 
our result implies that a compact isoperimetric region (when it exists) is bounded  by either a sphere or a torus. 
Under the weaker assumption that the scalar curvature is nonnegative we prove that if such a surface is connected then its genus is $\leq 3. $

In section 3, we particularize to the case when the ambient space is $\mathbb S^2\times\mathbb R$ and show that a closed and stable 
CMC surface is necessarily a rotational sphere.

In section 4, we address the problem in $\mathbb H^2\times\mathbb R.$ We obtain bounds on the genus of a stable closed CMC surface under
 various assumptions on the mean curvature. In particular we prove it is a rotational sphere if it has mean curvature 
$\geq 1/\sqrt{2}$ (Theorem \ref{h2timesR}). We believe this restriction on the mean curvature is unnecessary.

\section{STABILITY OF ROTATIONAL  CMC SPHERES IN  $\mathbb S^2\times \mathbb R$ AND  $\mathbb H^2\times \mathbb R$}
\medskip

The rotational CMC spheres in $\mathbb S^2\times \mathbb R$ and $\mathbb H^2\times\mathbb R$ have been well studied
 in the works of Hsiang and  Hsiang \cite{hsiang},  Pedrosa \cite{pedrosa} and  Abresch-Rosenberg \cite{abresch-rosenberg}.
In $\mathbb S^2\times\mathbb R,$  for any value $H>0$ there exists, up to ambient isometries, a unique rotational CMC
 sphere of mean curvature
$H.$ In $\mathbb H^2\times\mathbb R,$ they exist only for  $H>1/2.$ Moreover for any $H>1/2,$ there exists a unique, up to ambient isometries,
rotational CMC sphere of mean curvature $H$ in $\mathbb H^2\times\mathbb R.$  If then we fix the {\it center of gravity} of those surfaces to be at
some fixed point $\{p\} \times\{0\} \in \mathbb S^2\times\mathbb R$ (resp. in $\mathbb H^2\times\mathbb R$), then they form a
 smooth 1-parameter family $\Sigma_H$ parametrized by the mean curvature $H\in (0,+\infty)$ (resp. $H\in  (\frac{1}{2},+\infty)$). The surfaces $\Sigma_H$ are
embedded and invariant under the  reflection $\tau$ through the horizontal slice $\mathbb S^2\times\mathbb \{0\}$ (resp. $\mathbb H^2\times \{0\}$).
Furthermore the generating curve of any of these surfaces is a convex curve in a vertical geodesic plane.

We will need the following stability criterion due to  Koiso \cite{koiso}.

\begin{theorem}[\cite{koiso}, Theorem 1.3, statement (III-B)]\label{koiso}
Let $\Sigma$ be a smooth compact immersed and two-sided CMC surface in a Riemannian 3-manifold $M.$ Denote by $\lambda_1$ and $\lambda_2$ 
respectively the first and second eigenvalues of the Jacobi operator $L$ on $\Sigma.$ Suppose  the following conditions are 
verified:
\smallskip

 (i) $\lambda_1<0$ and $\lambda_2=0,$

(ii) $\int_{\Sigma} g \, dA =0$ for any eigenfunction $g$ associated to $\lambda_2.$

\smallskip

 In this case, denoting $E$ the eigenspace 
associated to $\lambda_2=0,$  there exists a uniquely determined smooth function $u\in E^{\perp}$ (the orthogonal to $E$ in the $L^2$ sense) which satisfies $Lu=1.$
Then: 
\smallskip

 (a) if  $\int_{\Sigma} u\,  dA\ge 0,$ then $\Sigma$ is stable,
 
 (b) if  $\int_{\Sigma} u\,  dA< 0,$ then $\Sigma$ is unstable.

\end{theorem}

Actually  Koiso states the theorem for compact CMC surfaces  with boundary in $\mathbb R^3$ but the
 result readily extends to  compact two-sided CMC surfaces with or without  boundary in  any ambient Riemannian 3-manifold. 
 The following result solves the problem of stability of rotational CMC spheres in $\mathbb S^2\times\mathbb R$ and
 $\mathbb H^2\times\mathbb R.$

\begin{theorem}\label{stablerotational}  Rotational CMC spheres in  $\mathbb H^2\times \mathbb R$ are stable.

\noindent There exists some value $H_0\approx 0.18$ such that rotational CMC spheres in  $\mathbb S^2\times \mathbb R$
are stable if  they have mean curvature $H\geq H_0$ and unstable if $0<H<H_0.$  Moreover a  finite union 
of horizontal slices in  $\mathbb S^2\times \mathbb R$ is stable. 

\end{theorem}

It is interesting to point out (cf. the proof below) that this critical value $H_0$ is the one where the area of the 
rotational CMC spheres reaches its maximum (it is also the value where the volume enclosed by these spheres is maximal).
 Pedrosa \cite{pedrosa}  has proved that, in $\mathbb S^2\times\mathbb R,$ the rotational CMC spheres of mean curvature 
$H\geq H_1\approx 0.33$ bound isoperimetric domains. Our result shows that the rotational CMC spheres are still stable  for 
$H_0\approx 0.18 \leq H\leq H_1\approx 0.33$ although they do not bound isoperimetric domains.

\begin{proof}

 We will apply Koiso's criterion (Theorem  \ref{koiso}).  We give details of the proof
 in the case of $\mathbb S^2\times \mathbb R,$ the case of $\mathbb H^2\times \mathbb R$ is quite similar. With the notations above, consider one of
the rotational CMC spheres $\Sigma_{H_1},$ for some $H_1>0$ ($\Sigma_{H_1}$  is rotationally invariant around the axis 
$\{p\} \times \mathbb R$ and has  center of gravity the point $\{p\} \times\{0\} \in \mathbb S^2\times\mathbb R$).  We will first check that conditions $(i)$ and $(ii)$ of Theorem \ref{koiso}
are verified.

$(i)$ First note that $0$ is an eigenvalue of the Jacobi operator $L.$  Indeed, if $X$ is any Killing field of $\mathbb S^2\times \mathbb R$
 then the mean curvature of the surface $\Sigma_{H_1}$ remains constant and equal to $H_1$ under the deformation  by flow of $X$ and so
  by (\ref{linearization})  the function $v= \langle X, N\rangle$ satisfies $Lv=0.$ 
   Recall that the isometry group of $\mathbb S^2\times\mathbb R$ has dimension 4 and so we have 3 independent Killing fields 
which do not vanish on $\Sigma.$ These fields induce on $\Sigma$ three linearly independent functions in the kernel of $L.$ 
These functions are indeed independent since otherwise the rotational surface $\Sigma_{H_1}$ would be invariant under a second
  1-parameter group of isometries of $\mathbb S^2\times\mathbb R,$ which is easily seen  to be not possible.  Therefore the dimension
 of $E=  ker
L$ is at least 3. As the first eigenspace of $L$ is one dimensional, we conclude that $\lambda_1<0.$

Given a point $q\in\mathbb S^2,$ we denote by $X_q$ the Killing field associated to the 1-parameter group of rotations around the 
vertical axis $\{q\}\times \mathbb R.$ Consider a vertical geodesic plane $\Pi$  containing the axis $\{p\}\times \mathbb R,$  that is $\Pi$ a
product $\Gamma\times\mathbb R \subset \mathbb H^2\times\mathbb R$ where $\Gamma$ is a geodesic of $\mathbb H^2$ containing $p.$ The plane $\Pi$ is a
plane of symmetry of $\Sigma_{H_1}$ and separates it into two components $\Sigma_{H_1}^1(\Pi)$ and 
$\Sigma_{H_1}^2(\Pi).$ Take a point $q\in\mathbb S^2,\, $ distinct from $p$ and from its antipodal point. Then
the function $\langle X_q,
N\rangle$ is a non zero function in the kernel of $L$ and vanishes precisely on the curve $\Pi\cap \Sigma_{H_1}.$
 In particular the first
eigenvalue of $L$ for the Dirichlet problem  on each of the components $\Sigma_{H_1}^1(\Pi)$  and
$\Sigma_{H_1}^2(\Pi)$ is equal to $0.$ 

In the same way, denoting by $t$ the last coordinate  in the product $\mathbb S^2\times\mathbb R,$
 the Killing field $\frac{\partial}{\partial t}$  induces a non  zero function $\langle \frac{\partial}{\partial t}, N\rangle$ in
 $E$ whose
zero set is the circle $\mathbb S^2\times\{0\}\cap  \Sigma_{H_1}.$
 This circle separates the surface into 2 components $\Sigma_{H_1}^+$ and
$\Sigma_{H_1}^-$ which are isometric by the reflection $\tau.$ Again we conclude
 that $0$ is the first eigenvalue of $L$ for the Dirichlet problem
on each of the components $\Sigma_{H_1}^+$ and $\Sigma_{H_1}^-.$

We now prove   that $\lambda_2=0.$ Assume by contradiction this is not the case, that is $\lambda_2<0.$
 Denote by $E_{\lambda_2}$ the eigenspace associated to $\lambda_2.$ The isometric involution $\tau$ on $\Sigma_{H_1}$
 (more precisely the
restriction of $\tau$ to   $\Sigma_{H_1}$) induces an involution on $E_{\lambda_2}$ and so this eigenspace splits
 into a direct sum of the
linear spaces of  invariant and anti-invariant  functions under the action of $\tau:$
$$ E_{\lambda_2}=\{v\in E_{\lambda_2}, \, v\circ \tau= v\} \oplus \{v\in E_{\lambda_2},\, v\circ \tau =-v\}.$$ 
Now any anti-invariant eigenfunction vanishes on the circle $\mathbb S^2\times\{0\}\cap \Sigma.$ Moreover by
 Courant's nodal domain theorem (\cite{cheng}) any non zero  eigenfunction in $E_{\lambda_2}$ has exactly 2 nodal domains.
Suppose then
there exists  a non zero anti-invariant eigenfunction in $E_{\lambda_2},$ its nodal domains are then precisely  
$\Sigma_{H_1}^+$ and $\Sigma_{H_1}^-.$
 It follows  that $\lambda_2$ is the first eigenvalue of $L$ for the Dirichlet problem on  $\Sigma_{H_1}^+$ (and
 on $\Sigma_{H_1}^-$) and this is a contradiction since we have seen that 0 is the first eigenvalue for this problem. So all functions in
$E_{\lambda_2}$ are invariant under the action of $\tau.$

Consider now as above  a vertical plane $\Pi=\Gamma\times\mathbb R$ containing the axis $\{p\}\times \mathbb R.$
 The surface $\Sigma_{H_1}$ is invariant under the symmetry $\sigma_{\Pi}$ through $\Pi$ and so as before $E_{\lambda_2}$ splits into a
direct sum of the spaces of invariant and anti-invariant functions, respectively, with respect to the induced action of 
$\sigma_{\Pi}.$ Again if we suppose there is a non zero anti-invariant eigenfunction it has to vanish on the curve 
$\Pi\cap\Sigma_{H_1}$
and as before this implies that $\lambda_2$ is an eigenvalue of $L$ for the Dirichlet problem on $\Sigma_H^1(\Pi)$ (and on
$\Sigma_H^2(\Pi)$) which is a contradiction. 
So all the functions in $E_{\lambda_2}$ are invariant with respect to reflection through $\Pi.$ As $\Pi$ was any vertical
 plane containing 
the axis $\{p\}\times \mathbb R,$ we conclude that each function in $E_{\lambda_2}$ is invariant under rotations through
 this axis. In particular the zero set of such an eigenfunction is also rotationally invariant. Courant's nodal domain
theorem then shows this zero set consists of exactly one  circle. As  any eigenfunction in $E_{\lambda_2}$ is invariant under the action of
$\tau$  this zero  set  is precisely the circle $\mathbb S^2\times\{0\}\cap \Sigma_{H_1}.$ This is again a contradiction. Consequently 
 $\lambda_2=0.$ 

$(ii)$   We can view the family $\Sigma_H, \, H\in(0,\infty)$ as
 given by a family of embeddings  $X_H$ of the surface $\Sigma_{H_1}.$ 
Call $Y$ the associated deformation vector field on $\Sigma_{H_1}$  and set $v= \langle Y, N\rangle.$
 Then by (\ref{linearization}) we have  
$Lv = 2.$ As $L$ is a self-adjoint operator, it follows that $\int_{\Sigma_{H_1}} g \, dA= 0,$  for any $g \in E= Ker L.$ Although we will not need  it, it is interesting to note that $E$ has dimension 3 and is generated by the functions induced by the Killing fields of the ambient space. Indeed by a result of Cheng \cite{cheng}, $E$ has dimension at most 3 and so coincides with the space induced by the Killing fields since  we have noted above that  this latter space has dimension 3. 

Write  $\frac{1}{2} v  = u + h$ with $u\in E^{\perp}$ and $h\in E.$ Then $Lu = 1$ and $\int_{\Sigma_{H_1}} u \, dA=\frac{1}{2}  \int_{\Sigma_{H_1}} v\, dA .$ So to conclude the proof, it remains to study the sign of $\int_{\Sigma_{H_1}}v\,  dA.$

 Call $A(H)$ the area
of  $\Sigma_H.$ Then by the first variation formula for the area: $\frac{dA}{dH}(H_1)= -2H_1\int_{\Sigma_{H_1}} v\, dA.$ 
Therefore the sign of $\int_{\Sigma_{H_1}} v \, dA$
is the opposite of that of $\frac{dA}{dH}(H_1).$ Now we have an explicit expression for $A(H),$ derived by  Pedrosa \cite{pedrosa}.
 Namely:
$$A(H)= 8\pi\left[    \frac{1}{4H^2+1}+\frac{4H^2}{(4H^2+1)^{3/2}}\tanh^{-1} \frac{1}{\sqrt{4H^2 +1}}\right]$$
(Pedrosa takes as definiton of the mean curvature the trace of the second fundamental form while we take half of it).
A study of the derivative function $\frac{dA}{dH}$  shows it has a unique zero on  $(0,\infty)$ at some value $H_0\approx 0.18$ and
 is positive on $(0,H_0)$ and negative on $(H_0,\infty).$ Therefore the surfaces $\Sigma_H$ verify  condition $(a)$ of Theorem 
(\ref{koiso}) for $H\geq H_0$ and condition $(b)$ for $0<H<H_0.$ This concludes the proof for the case of
 $\mathbb S^2\times \mathbb R.$

In $\mathbb H^2\times\mathbb R$ the expression for the area $A(H)$ was derived by  Hsiang and Hsiang \cite{hsiang}. In this case, one has:
$$A(H)= 8\pi\left[\frac{1}{4H^2 -1}+\frac{4H^2}{(4H^2-1)^{3/2}} \tan^{-1}\frac{1}{\sqrt{4H^2 -1}}\right].$$
It can be checked that the derivative function $\frac{dA}{dH}$  is always negative. The conclusion follows
 from case $(a)$ of Theorem \ref{koiso}. 

The last statement is quite easy and in fact a finite union of horizontal slices in 
$\mathbb S^2\times\mathbb R$ is  stable in a strong sense. Indeed, as  a horizontal slice is totally geodesic and the Ricci curvature evaluated
on its normal vanishes, its Jacobi operator is  just the Laplacian on the sphere. Therefore the stability condition (\ref{stability}) is
satisfied for any smooth function defined on a finite union of horizontal slices without any condition on its integral.
\end{proof}

\begin{remark} The stability of rotational CMC spheres in $\mathbb H^2\times\mathbb R$ can be seen in an alternate way.
Actually  they bound isoperimetric domains \cite{hsiang}. We think, however, that  the direct argument we gave
here is interesting on its own. 
\end{remark}

\section{STABLE  CMC SURFACES IN SIMPLY CONNECTED CONFROMALLY FLAT  MANIFOLDS}

 We obtain in this section some control on the topology of stable compact and orientable CMC surfaces in 
simply connected conformally flat manifolds. Recall that a Riemannian
manifold is said to be  {\it conformally flat\rm}
 if each of its points
admits an open neighborhood which is conformally diffeormorphic to an open set of the Euclidean space. 
We will need a first result.

\begin{proposition} \label{conformal}Let  $M$ be  a simply connected conformally flat Riemannian 3-manifold.
 Denote by $K_s$ its sectional curvature. Then for any compact 
orientable surface $\Sigma$ without boundary immersed in $M$ we have:
\begin{equation}\label{willmore}
\int_{\Sigma} \{H^2(x) + K_s(T_x\Sigma)\} dA(x) \geq 4\pi,
\end{equation}
where $H$ is the mean curvature of $\Sigma,$ and   equality holds if and only if $\Sigma$ is a totally umbilic sphere.

Furthermore if $\Sigma$ is not embedded then:
\begin{equation}\label{li-yau}
\int_{\Sigma} \{H^2(x) + K_s(T_x\Sigma)\} dA(x) \geq 8\pi.
\end{equation}

\end{proposition}

\begin{proof}
It is known that   the integral $\int_{\Sigma} \{H^2(x) + K_s(T_x\Sigma)\} dA(x)$
is invariant under conformal changes of the metric on $M.$ It is also known that umbilicity of a surface in a 3-manifold
 is preserved by conformal changes of the metric on the ambient manifold (cf. for instance \cite{spivak}). Moreover
 for a simply connected
conformally flat Riemannian 3-manifold the {\it developing map \rm} $\text{dev}:  M \longrightarrow \mathbb S^3$ is a well defined
conformal immersion.  It is therefore enough to check the theorem when the ambient manifold is the canonical 3-sphere $\mathbb S^3.$ The first
statement  is then due to Willmore \cite{willmore}, and the second one is a result of Li and Yau  \cite{L-Y}. 
\end{proof}
\begin{remark}\label{remark}
 Proposition \ref{conformal} applies in particular to the product spaces $\mathbb S^2\times \mathbb R$ and $\mathbb H^2\times\mathbb R$
 as they are conformally flat. More precisely $\mathbb S^2\times \mathbb R$ is conformal to $\mathbb R^3\setminus\{(0,0,0)\}$
and $\mathbb H^2\times \mathbb (0,\pi)$ is conformal to $\mathbb H^3$ (cf. \cite{souam-toubiana}).
 The totally umbilic surfaces in these spaces and more generally in homogeneous 3-manifolds were classified recently
 by Toubiana and the
author in \cite{souam-toubiana2}.
\end{remark}
Our main result in this section is:

\begin{theorem} \label{thm:conformal}
Let $M$ be a simply connected conformally flat Riemannian three manifold. 
Let $\Sigma$ be a compact orientable and  stable 
CMC surface  in $M$ with unit normal $N.$

(i) Suppose the Ricci curvature of $M$ is nonnegative, then either $\Sigma$ is a
 finite union of totally geodesic surfaces and Ric($N)=0$  or $\Sigma$ is a  sphere or an embedded torus.
 
(ii) Suppose the scalar curvature of $M$ is nonnegative and that $\Sigma$ is  connected. Then $\Sigma$ has
genus $g\leq 3.$ Moreover,  if $g=3,$  then  $\Sigma$ is  embedded and if
$g=2$ and $\Sigma$ is not embedded then it is minimal  and the scalar curvature of $M$ vanishes on $\Sigma.$  
\end{theorem}

\begin{proof} (i) We first note that if $\Sigma$ is not connected, then we can take a locally constant function with
 vanishing integral and easily contradict stability if the function $ |\sigma|^2+{\text Ric}(N)$ is not identically zero on $\Sigma.$

Assume now that $\Sigma$ is a connected surface of genus $g.$
By Brill-Noether's theory (cf. \cite{G-H})
there exists a  non constant conformal (i.e holomorphic)  map $\phi: \Sigma \longrightarrow \mathbb S^2 \subset \mathbb R^3$ 
 whose degree satisfies:
\begin{equation}\label{degree}
 \text{deg} (\phi) \leq 1+\left[\frac{g+1}{2}\right] 
 \end{equation}
where $[x]$ denotes the integer part of $x.$ 
Using an extended version of a result of Hersch \cite{hersch, L-Y}) we may assume after composing 
$\phi$ with a conformal diffeomorphism of $\mathbb S^2$ that its coordinate functions 
satisfy:
\begin{equation}\label{meanvalue}
 \int_{\Sigma } \phi_i \,dA =0,\,\,\, i=1,2,3.
 \end{equation}
We may then use these coordinate functions as test functions for the stability condition and obtain:

\begin{equation}\label{holom}
  0\leq Q(\phi_i,\phi_i)= \int_{\Sigma} |\nabla \phi_i|^2 - (|\sigma|^2+\text{Ric}(N))\phi_i^2 ,\quad\quad i=1,2,3.
  \end{equation}
Summing up these inequalities and taking  into account that $\phi$ takes its values in  the unit sphere, we get:

\begin{equation} \label{eq:holo}
0 \leq  \int_{\Sigma} |\nabla \phi|^2 - (|\sigma|^2+\text{Ric}(N)) 
\end{equation}
Let $K_s$ and $K_{\Sigma}$ denote, respectively, the sectional curvature of $M$ evaluated on the tangent plane to $\Sigma$  and the intrinsic curvature 
of $\Sigma$. By Gauss equation we have: $|\sigma|^2 = 4H^2 +2K_s -2 K_{\Sigma}.$ Using Gauss-Bonnet theorem and 
the fact that for a conformal map $\int_{\Sigma} |\nabla \phi|^2 = \int_{\Sigma}  2\text{Jacobian}(\phi)= 8\pi \text{deg}(\phi),$ we transform inequality (\ref {eq:holo}) into:
\begin{equation}\label{eq:holo1}
\int_{\Sigma} 4H^2 +2 K_s + \text{Ric(N)} \leq 8\pi\left(2-g+\left[\frac{g+1}{2}\right]\right ).
\end{equation}
Using inequality (\ref{willmore}) of Proposition \ref{conformal}, we get:
\begin{equation}\label{g=2 or 3}
\int_{\Sigma} 2H^2 + \text{Ric}(N) \leq 8\pi\left(1-g+\left[\frac{g+1}{2}\right]\right).
\end{equation}
As by hypothesis  the Ricci curvature is nonnegative, this implies $g\leq 3.$ 

If $g=2$ or $3,$ then  the right-hand side in (\ref{g=2 or 3}) is equal to zero and since the Ricci curvature is nonnegative
 there must be equality in (\ref{g=2 or 3}). Furthermore all  the intermediate  inequalities have to be  equalities and in particular we
have the  equality in (\ref{willmore}) which by Proposition  \ref{conformal} implies $\Sigma$ is a (totally umbilic) sphere and that is a
contradiction. 

Suppose now $\Sigma$ is a torus and assume it is not embedded. Then  taking into account  this time inequality (\ref{li-yau}) 
in Proposition \ref{conformal} we arrive at the inequality:
\begin{equation} \label{torus}
\int_{\Sigma} 2H^2 +\text{Ric}(N) \leq 0.
\end{equation}
Again as the Ricci curvature is nonnegative there must be   equality   in  (\ref{torus}) and  so all the intermediate
  inequalities have to be equalities too. In particular  we have equalities in (\ref{holom}). So the holomorphic map $\phi$ satisfies
$Q(\phi_i,\phi_i)=0,$ for $ i=1,2,3.$ As $\Sigma$ is stable, for any $v\in\mathcal{C}^\infty (\Sigma)$ satisfying $\int_{\Sigma} v=0$
and any $t\in\mathbb R$ we have $Q(\phi_i +tv,\phi_i+tv) \geq 0$ and so $Q(\phi_i,v)=0.$ It follows  that each of the  functions $\phi_i$
 satisfies the equation 
$$\Delta \phi_i + |\sigma|^2 \phi_i =c_i,$$
for some real constant $c_i,\, i=1,2,3.$ So $\phi$ satisfies an equation of the type:
\begin{equation}\label{phi1}
\Delta \phi + |\sigma|^2 \phi=c
\end{equation}
with $c$ a constant vector in $\mathbb R^3.$ 

Since $\phi: \Sigma \longrightarrow \mathbb S^2$ is holomorphic it is harmonic and therefore satisfies the equation:
\begin{equation}\label{phi2}
\Delta \phi + |\nabla \phi|^2 \phi=0.
\end{equation}
As $\phi$ takes its values in the sphere $\mathbb S^2$ and is non constant, it follows easily from (\ref{phi1})
 and (\ref{phi2}) that necessarily $c=0$ and $|\sigma|^2=|\nabla\phi|^2.$ So the Jacobi operator of
$\Sigma$ writes as $L= \Delta +|\nabla\phi|^2$ and the stability assumption implies that $L$ has only one negative eigenvalue.
 Otherwise said the holomorphic map $\phi$
 has index one and such maps do not exist on tori (cf.  \cite{ros2}). So if it is a torus, $\Sigma$ has to be embedded. 
This completes the proof of the first part of the theorem. 

(ii) To prove the second part of the theorem, we proceed in the same way. Taking into account that Ric$(N)= S-K_s,$
where $S$ denotes the scalar curvature of $M,$ we rewrite (\ref {eq:holo1}) as follows:
\begin{equation}\label{eq:holo3}
 \int_{\Sigma} 4H^2 + K_s + S \leq 8\pi\left(2-g+\left[\frac{g+1}{2}\right]\right ).
 \end{equation}
Using again inequality (\ref{willmore}) we obtain:
$$\int_{\Sigma} 3H^2+S\leq 4\pi \left( 3-2g+2\left[\frac{g+1}{2}\right]\right ).$$
As $S\geq 0$ by hypothesis, this implies $g\leq 3.$ Furthermore, suppose $\Sigma$ is not embedded then using (\ref{li-yau})
in (\ref{eq:holo3}) we get:
\begin{equation}\label{eq:holo4}
 \int_{\Sigma} 3H^2 + S \leq 8\pi\left(1-g+\left[\frac{g+1}{2}\right]\right ).
 \end{equation}
If then $g=2$ or $3,$ the right-hand side in (\ref{eq:holo4}) is zero and so $H=0$ and $S=0$ on $\Sigma.$ Moreover 
the equality is achieved in (\ref{eq:holo4}) and therefore  also in the intermediate inequalities. As above we conclude
 that the Jacobi operator writes as $L=\Delta+ |\nabla \phi|^2$ and thus the holomorphic map $\phi:\Sigma\longrightarrow \mathbb S^2$
has index one. We have deg$(\phi)=3$ when $g=3.$  But on a surface of genus 3, a holomorphic map of index one necessarily
 has degree 2 (cf.  \cite{ros2} or \cite{ritore}). So if $g=3$ the surface has to be embedded.
\end{proof}

In particular in case (i) if we assume the Ricci curvature is positive,  a compact stable CMC surface
can be either a sphere or a torus.  As a consequence a compact  isoperimetric region in a simply connected locally conformally 
flat 3-manifold with positive Ricci curvature is bounded by either a sphere or a torus.

\section{ STABLE CMC SURFACES IN $\mathbb S^2\times \mathbb R$ }

Theorem \ref{thm:conformal} applies in particular to $\mathbb S^2\times \mathbb R$ (cf. Remark  \ref{remark}).
 However in this case using the symmetries of this space we can give a complete solution to the stability problem. 

\begin{theorem}\label{s2timesR}
Let $\Sigma$ be a compact stable immersed CMC surface without boundary in $\mathbb S^2\times \mathbb R.$ Then $\Sigma$ is either a finite union
 of horizontal slices or a rotational sphere of mean curvature 
$H\geq H_0\approx 0.18.$

\end{theorem}

\begin{proof} It follows from the maximum  principle that the horizontal slices $\mathbb S^2\times \{t\},\, t\in\mathbb R,$ 
are the only compact minimal surfaces without boundary in $\mathbb S^2\times\mathbb R.$ In particular all compact CMC surfaces in $\mathbb S^2\times
\mathbb R$ are orientable.  By  Theorem \ref{thm:conformal} (and Remark \ref{remark}),  if $\Sigma$ is not connected it is a finite union of horizontal slices and if  connected it is either a sphere or an embedded torus. We
will discard the case of a torus. Assume such a stable torus $\Sigma$ exists. As it has to be embedded, a simple application of Alexandrov
reflection technique through the  reflections $(p,t)\in \mathbb S^2\times \mathbb R \to (p, 2t_0-t),\, t_0\in \mathbb R,$ shows that $\Sigma$ is
invariant under one such reflection.   Call $\tau$ this reflection and let  $\Sigma_+=\Sigma\cap\big(\mathbb S^2\times [t_0,+\infty)\big).$
Note that $\Sigma_+$ is topologically an annulus.

The argument to follow uses the symmetries of the ambient space and is based on Courant's  nodal domains theorem. A 
similar idea was used to study volume-preserving stability in some variational problems \cite{ros-vergasta, 
 ros-souam, souam, hutchings et al}.  Let 
$q=(p_1,t_1)$ be a point on $\Sigma$ where the height function reaches its maximum. Thus $q\in \Sigma^+.$ Denote by $N$ a global unit
normal on $\Sigma.$ Consider the 1-parameter  family of rotations $f_{\theta}$ of $\mathbb S^2\times \mathbb R$ around the axis
$\{p_1\}\times\mathbb R.$ It follows from (\ref{linearization}) that the associated infinitesimal  rotation $u:=\langle \frac{df}{d\theta},
N\rangle$ is a Jacobi function, i.e it satisfies the equation:
\begin{equation}\label{jacobi}
\Delta u + (|\sigma|^2+\text{Ric}(N)) u = 0.
\end{equation}

Also it is straightforward to check that 
\begin{equation}\label{critical} u(q)=0\quad \text{and} \quad \nabla u(q)=0.
\end{equation}
We will show that $u\equiv 0$ on $\Sigma.$ Assume by contradiction this is not the case. 
Then as $u$ satisfies (\ref{jacobi}), by a result of Cheng \cite{cheng} its zero set $\{u=0\}$ is a (maybe disconnected) graph with
vertices the points where $\nabla u$ vanishes too. Moreover any sufficiently small neighborhood of  a vertex  is divided
 by the zero set of
$u$  like a pie, into at least 4  regions, such that in any two adjacent regions $u$ has opposite signs (\cite{cheng}). By
(\ref{critical}), the point $q\in \Sigma_+$ is one such  a vertex.  It is a consequence of the Jordan curve theorem that $\Sigma_+
\setminus\{u=0\}$ has at least three connected components. If now $\Omega \subset \Sigma_+$ is one such component then either $\Omega$
 and
$\tau(\Omega)$ are two distinct connected components of $\Sigma\setminus\{u=0\}$ or $\Omega\cup\tau(\Omega)$ is a connected component 
of 
 $\Sigma\setminus\{u=0\}.$ Therefore the number of connected components of $\Sigma\setminus\{u=0\}$ is at least as big as the number
 of connected components of $\Sigma_+ \setminus\{u=0\}.$ In particular $\Sigma \setminus\{u=0\}$ has at least three connected components. 
 Denote by $\lambda_1< \lambda_2$ respectively the first and second eigenvalues of the operator $L$ on $\Sigma.$
 It is well known that $\lambda_1$ is simple and if  $u_1$ is a (non zero)
 associated eigenfunction, then  $u_1$ does not vanish on $\Sigma.$  Also by Courant nodal domain theorem (cf.\cite{cheng}),
if  
$u_2$ is a non zero  eigenfunction associated to $\lambda_2,$
 then the set $\Sigma\setminus\{u_2=0\}$ has exactly two connected components. As $u$ is an eigenfunction associated
 to the eigenvalue zero of $L,$ by  (\ref{jacobi}), it follows from the previous considerations that $\lambda_2 <0.$
 Let $a=-\int u_2/\int
u_1,$
 and set $v:= a u_1+ u_2.$  Then  $\int_{\Sigma} v=0$ and:
 $$Q(v,v)= a^2 Q(u_1,u_1)+ Q(u_2,u_2)= a^2\lambda_1\int_{\Sigma} u_1^2 +\lambda_2 \int_{\Sigma} u_2^2 <0$$
 and this contradicts stability.

Therefore $u\equiv 0,$ which means that $\Sigma$ is a surface of revolution  around the axis $\{p_0\}\times \mathbb R.$
As this axis intersetcs $\Sigma,$ this excludes the case of a torus. 

The only left possibility is therefore that $\Sigma$ is a sphere. To conclude that $\Sigma$ is a sphere of revolution
 we can use the Abresch-Rosenberg theorem \cite{abresch-rosenberg} which asserts that an immersed CMC sphere in 
$\mathbb S^2\times\mathbb R$ or in $\mathbb H^2\times \mathbb R$ is rotationally invariant. Alternatively we can also
repeat  the previous argument
 to see that it is
a sphere of revolution around the axis passing through a point where the height function reaches its maximum,
the existence of the symmetry
is not needed when $\Sigma$ is a sphere.  We know from  Theorem
\ref{stablerotational} that the mean curvature of a stable CMC sphere is $\geq H_0\approx 0.18.$
\end{proof}

In particular the Theorem says there is no compact stable CMC
 surface in $\mathbb S^2\times \mathbb R$ of mean curvature $0<H<H_0\approx 0.18.$

\section{ STABLE CMC SURFACES IN $\mathbb H^2\times \mathbb R$}

\medskip

In this section we study compact stable CMC surfaces in $\mathbb H^2\times\mathbb R.$
 Nelli and Rosenberg \cite{nelli-rosenberg}
have proven that the surface has genus at most 3
 if its mean curvature is $>1/\sqrt 2.$ In view of the result obtained in the
case of $\mathbb S^2\times \mathbb R,$ it is reasonable to conjecture that only rotational CMC spheres are stable. Our next result gives a 
partial answer:

\begin{theorem} \label{h2timesR} Let $\Sigma$ be a compact immersed and connected surface without boundary of constant mean curvature $H$ and genus $g$
in $\mathbb H^2\times
\mathbb R.$ Assume $\Sigma$ is stable. 

\noindent If 
$H\geq 1/\sqrt{2},$ then $\Sigma$ is a rotational sphere.

 \noindent If $1/\sqrt {3} <H<1/\sqrt{2},$ then $g\leq 1.$
 
\noindent If $H=1/\sqrt{3}$ then $g\leq 2.$
\end{theorem}
\medskip
\begin{remark} Nelli and Rosenberg \cite{nelli-rosenberg} have proved that an immersed 
compact CMC surface without boundary in $\mathbb
H^2\times \mathbb R$ always has mean curvature 
$H> 1/2.$ 
\end{remark}
\begin{proof}
First note that by the maximum principle there are no compact minimal surfaces without boundary in $\mathbb H^2\times\mathbb R$
and so in particular $\Sigma$ is orientable. We proceed as in the case of $\mathbb S^2\times \mathbb R.$ So we have
 a non constant holomorphic map $\phi: \Sigma\longrightarrow \mathbb S^2\subset\mathbb R^3$ satisfying conditions (\ref{degree}) and
(\ref{meanvalue}) and this leads to the inequality:
\begin{equation}\label{eq:holo2}
\int_{\Sigma} 4H^2 +2 K_s + \text{Ric}(N) \leq 8\pi\left(2-g+\left[\frac{g+1}{2}\right] \right).
\end{equation}
The Alexandrov reflection technique through planes of the form $\Gamma \times \mathbb R,$  where $\Gamma$ is any geodesic in 
$\mathbb H^2,$ shows that the CMC spheres of revolution  are the only embedded compact CMC surfaces without boundary in 
$\mathbb H^2\times\mathbb R.$ So keeping aside the genus zero case, we have that   $\Sigma$ is not embedded. Using
 then the inequality (\ref{li-yau}) in (\ref{eq:holo2}), we get:
\begin{equation}
\int_{\Sigma} 2H^2+\text{Ric}(N) \leq 8\pi\left(-g+\left[\frac{g+1}{2}\right]\right)
\end{equation}
In $\mathbb H^2\times\mathbb R$ we have  $-1\leq$  Ric$(V)\leq 0$ for any unit vector $V$ and Ric$(V)=-1$
if and only if $V$ is horizontal. So:
\begin{equation}
\int_{\Sigma} 2H^2-1\leq \int_{\Sigma} 2H^2+\text{Ric}(N) \leq 8\pi\left(-g+\left[\frac{g+1}{2}\right]\right) \leq 0
\end{equation}
Consequently if $H> 1/\sqrt{2}$ we reach a contradiction. 

If now $H=1/\sqrt{2}$
then necessarily Ric$(N)=-1$ and so   $N$ is always horizontal
which is absurd as $\Sigma$ is compact without boundary. So if $H\geq 1/\sqrt{2}$ the surface has to be
 a sphere and we can conclude as in the proof of Theorem \ref{s2timesR} that it is a sphere of revolution.

To obtain the second part of the theorem, note that since the scalar curvature of $\mathbb H^2\times\mathbb R$
is identically $-1$  we have  Ric$(N)=-1-K_s.$ So we can rewrite (\ref{eq:holo2}) as follows:
\begin{equation}
\int_{\Sigma} 4H^2 +K_s -1 \leq 8\pi\left(2-g+\left[\frac{g+1}{2}\right] \right).
\end{equation}
Again, keeping aside the genus zero case, we know the surface is non embedded and so we can use inequality (\ref{li-yau}) to obtain:
\begin{equation}\label{3h2}
\int_{\Sigma} 3H^2-1 \leq 8\pi\left(1-g+\left[\frac{g+1}{2}\right] \right).
\end{equation}
This shows $g\leq 1$ when $H>1/\sqrt{3}.$

Suppose now that $H=1/\sqrt{3},$ then by (\ref{3h2}) we have $g\leq 3.$ Moreover if $g=2$ or $3,$ we have equality in 
(\ref{3h2}) and in all the intermediate inequalities and we conclude as in the proof of Theorem \ref{s2timesR}
that the Jacobi operator writes as $L= \Delta + |\nabla\phi|^2$ and so the holomorphic map $\phi$ has index one. Note that
deg$(\phi)=3$ when $g=3.$ However  on a surface of genus three a holomorphic map of index one necessarily has degree two (cf.\cite{ros2}
or \cite{ritore}). This excludes the case $g=3.$
\end{proof}

\bibliographystyle{alpha}

\begin{thebibliography}{06}



\bibitem{abresch-rosenberg} U. Abresch \& H. Rosenberg.: 
{\em  A Hopf differential
for constant mean curvature surfaces in $\mathbb S^2\times \mathbb R$ and $\mathbb H^ 2\times \mathbb R$. } Acta Math.  193 
(2004),  no. 2, 141--174.


\bibitem{barbosa-do carmo} J-L. Barbosa \& M. do Carmo.:
{\em  Stability of
hypersurfaces with constant mean curvature.}  Math. Z. 
185  (1984),  no. 3, 339--353.

\bibitem{barbosa et al} J-L. Barbosa; M.  do Carmo \& J. Eschenburg.:
{\em 
 Stability of hypersurfaces of constant mean
curvature in Riemannian manifolds. } Math. Z.  197 
(1988),  no. 1, 123--138.


\bibitem{cheng} S.Y. Cheng.:
 {\em  Eigenfunctions and nodal sets.}
Comment. Math. Helv.  51  (1976), no. 1, 43--55.


\bibitem{G-H} P. Griffiths \& J. Harris 
{\em  Principles of algebraic geometry. Pure and Applied Mathematics.}
Wiley-Interscience, John Wiley \& Sons], New York, 1978.

\bibitem{hersch} J. Hersch.:
{\em 
Quatre propri\'et\'es isop\'erim\'etriques de membranes
sph\'eriques homog\`enes. }
C. R. Acad. Sci. Paris SŽr. A-B 270 (1970),
A1645--A1648.

\bibitem{hsiang} W-T. Hsiang \& W-Y. Hsiang.:
{\em  On the uniqueness of isoperimetric solutions and imbedded soap bubbles in noncompact symmetric spaces. I.}  Invent. Math.  98  (1989),  no. 1, 39--58. 

\bibitem{hutchings et al} M. Hutchings; F. Morgan; M. Ritor\'e \& A. Ros.: 
{\em  Proof of the double bubble conjecture.}
Ann. of Math. (2)  155  (2002),  no. 2, 459--489



\bibitem{koiso} M. Koiso.:
{\em  Deformation and stability of surfaces
with constant mean curvature.}  Tohoku Math. J. (2)  54
 (2002),  no. 1, 145--159.


\bibitem{L-Y} P. Li \& S-T. Yau.:
{\em  A new conformal invariant
and its applications to the Willmore conjecture and
the first eigenvalue of compact surfaces. } Invent.
Math.  69  (1982), no. 2, 269--291

\bibitem{nelli-rosenberg} B. Nelli \& H. Rosenberg.:
{\em  Global properties of constant mean curvature surfaces 
in $\mathbb H^ 2\times\mathbb R$.}  Pacific J. Math.  226  (2006),  no. 1, 137--152. 

\bibitem{pedrosa} R. Pedrosa.: 
{\em  The isoperimetric problem in spherical cylinders.}  Ann. Global Anal. Geom.  26  (2004),  no. 4, 333--354. 



\bibitem{ritore} M. Ritor\'e.: 
{\em Index one minimal surfaces in flat
three space forms.}  Indiana Univ. Math. J.  46 
(1997),  no. 4, 1137--1153.

\bibitem{ros1} A. Ros.: 
{\em The isoperimetric problem.} Global theory
of minimal surfaces,  175--209, Clay Math. Proc., 2,
Amer. Math. Soc., Providence, RI, 2005.


\bibitem{ros2} A. Ros.:
{\em  One-sided complete stable minimal
surfaces. } J. Differential Geom.  74  (2006),  no. 1,
69--92. 


\bibitem{ros3} A.  Ros.: {\em Stable periodic constant mean curvature surfaces and mesoscopic phase separation.}  Interfaces and Free Boundaries (to appear).

\bibitem{ros-souam} A. Ros \& R. Souam.:
{\em  On stability of capillary
surfaces in a ball.}  Pacific J. Math.  178  (1997), 
no. 2, 345--361.

\bibitem{ros-vergasta} A. Ros \& E. Vergasta.: {\em Stability for hypersurfaces of constant mean curvature with free boundary.}  Geom. Dedicata  56  (1995),  no. 1, 19--33.

\bibitem{ross} M. Ros.:  {\em Schwarz' $P$ and $D$ surfaces are stable.}  Differential Geom. Appl.  2  (1992),  no. 2, 179--195.

\bibitem{souam} R. Souam.:  {\em  On stability of stationary hypersurfaces
for the partitioning problem for balls in space forms.}
 Math. Z.  224  (1997),  no. 2, 195--208. 

\bibitem{souam-toubiana}R. Souam \& E. Toubiana.: {\em On the classification and regularity of umbilic surfaces in homogeneous 3-manifolds.} Matem\`atica Contempor\^anea (to appear).

\bibitem{souam-toubiana2}R. Souam \& E. Toubiana.: {\em Totally umbilic surfaces in homogeneous 3-manifolds.} Preprint (2006).

\bibitem{spivak} M. Spivak. {\em A comprehensive introduction to
    differential geometry}, Vol. 4, Boston, Publish or Perish, 
1970.

\bibitem{willmore} T.J. Willmore.: {\em  Note on embedded surfaces.}  An. \c Sti. Univ. "Al. I. Cuza" Ia\c si Sec\c t. I a Mat. (N.S.)  11B  1965 493--496. 
\end{thebibliography}

\enddocument